\documentclass[draft]{amsart}
\usepackage{amssymb} 
\usepackage[utf8]{inputenc}
\newtheorem{thm}{Theorem}[section]
\newtheorem{cor}[thm]{Corollary}

\newtheorem{pro}[thm]{Proposition}

\newtheorem{exa}[thm]{Example}

\newtheorem{rem}[thm]{Remark}
\numberwithin{equation}{section}

\begin{document}

\title[$C_{BSE}(\Delta(A))$ AS A DUAL BANACH ALGEBRA]
{$C_{BSE}(\Delta(A))$ AS A DUAL BANACH ALGEBRA}%

\author{Fatemeh Abtahi, Ali Rejali and Farshad Sayaf}%


%




\keywords{BSE functions, Dual Banach algebra, Lipschitz algebra.}

\subjclass[2000]{46A03, 46A04}


\begin{abstract}
Let $(X,d)$ be a metric space and $A$ be a commutative Banach algebra such that $\Delta(A)$
be nonempty. In this paper, we provide necessary and sufficient conditions, under which
$C_{BSE}(\Delta({A}))$ is a Dual Banach algebra. Moreover, we provide some conditions 
for that $Lip_\alpha(X,A)$ is a Dual Banach algebra. 
\end{abstract}

\maketitle

\section{\bf Introduction and Preliminaries}

Let $(A,\|\cdot\|)$ be a commutative Banach algebra and ${A}'$ and ${A}''$ be the dual and second dual Banach spaces, respectively. Let $a\in A$, $f\in {A}'$ and $\Phi,\Psi\in {A}''$. Then $f\cdot a$ and $a\cdot f$ are defined as $f\cdot
a(x)=f(ax)$ and $a\cdot f(x)=f(xa)$, for all $x\in A$,
making ${A}'$ an $A-$bimodule. Moreover for all
$f\in {A}'$ and $\Phi\in {A}''$, we define
$\Phi\cdot f$ and $f\cdot\Phi$ as the elements of ${A}'$ by
$$
\langle\Phi\cdot f,a\rangle=\langle\Phi,f\cdot a\rangle\;\;and\;\;
\langle f\cdot \Phi,a\rangle=\langle\Phi,a\cdot
f\rangle\;\;\;\;\;\;\;\;\;\;\;(a\in A).
$$
This defines two Arens products $\square$ and $\lozenge$ on
${A}''$ as
$$
\langle\Phi\square\Psi,f\rangle=\langle\Phi,\Psi\cdot
f\rangle\;\;and\;\;
\langle\Phi\lozenge\Psi,f\rangle=\langle\Psi,f\cdot \Phi\rangle,
$$
making ${A}''$ a Banach algebra with each. The products
$\square$ and $\lozenge$ are called respectively, the first and
second Arens products on ${A}''$. Note that $A$
is embedded in its second dual via the identification
$$
\langle a,f\rangle=\langle f,a\rangle\;\;\;\;\;\;\;\;(f\in
{A}');
$$
see \cite{D}, for more information. The space, consisting of all nonzero multiplicative linear functionals on $A$ will be denote by  $\Delta({A})$ and is called the Gelfand space
of $A$. A bounded continuous function $\sigma$ on $\Delta({A})$ is called a
BSE-function if there exists a constant $M>0$ such that for every finite number of
$\varphi_1,\cdots,\varphi_n$ in $\Delta({A})$ and the same number of $c_1,\cdots,c_n$ in $\Bbb C$, the inequality
$$\left\vert \displaystyle\sum_{j=1}^nc_j\sigma(\varphi_j)\right\vert\leq M\;\left\Vert
\displaystyle\sum_{j=1}^nc_j\varphi_j\right\Vert_{{A}^*}$$ holds.  The BSE-norm of
$\sigma$ is the infimum of all such $M$ and is denoted by $\Vert\sigma\Vert_{BSE}$.
The set of all BSE-functions  will be denoted by $C_{BSE}(\Delta({A}))$.
It has been shown that  under the norm $\Vert\cdot\Vert_{BSE}$,
$C_{BSE}(\Delta({A}))$ is a commutative and semisimple Banach algebra, 
embedded in $C_b(\Delta({A}))$ as a subalgebra. We refer to \cite{Tak1}, for a
full information about the history of BSE functions and BSE algebras.

Now we provide some information about Lipschitz algebras;
see \cite{EM}, \cite{HNS} and \cite{J}, for more information.
Let $(X,d)$ be a metric space and $A$ be a Banach algebra. For $f:X\rightarrow A$ let
$$
\|f\|_{\infty,E}=\sup_{x\in X}\|f(x)\|.
$$
Then $f$ is called bounded if $\|f\|_\infty<\infty$.
Now let
$$
\rho_{d,E}(f)=\sup_{x\neq y}\displaystyle\frac{\|f(x)-f(y)\|_E}{d(x,y)}.
$$
Following \cite{EM} and \cite{HNS}, the vector-valued Lipschitz space
$Lip_d(X,A)$, is the space consisting of all bounded
maps $f:
X\rightarrow A$ with $\rho_{d,E}(f)<\infty$. Furthermore,
$Lip_d(X,A)$, equipped with the norm $$\|f\|=\max\{\rho_{d,E}(f),\|f\|_{\infty,E}\}$$ and
pointwise product, is a Banach algebra.
Moreover,  we denote by $lip_d(X,A)$, the subspace of
$Lip_d(X,A)$, consisting of all functions $f$ such
that
$$\lim_{d(x,y)\rightarrow
0}\displaystyle\frac{\|f(x)-f(y)\|}{d(x,y)}=0.$$
Then $lip_d(X,A)$ is a Banach algebra with
pointwise product. In the case where $A$ is the field of complex numbers $\Bbb C$, to
simplify the notation, we will write $Lip_d X$ and $lip_d X$, rather than $Lip_d(K,\Bbb C)$ and $lip_d(X,\Bbb C)$, respectively. 

Recently,  some important algebraic properties of Lipschitz algebras have been investigated; see
\cite{AAR}, \cite{ABR} and \cite{AKT}. Here, we study another feacure of this algebras. 
We first recall the notion of Dual Banach algebras.

Following \cite{R}, a dual Banach algebra is a pair $(A,A_*)$ such that:
\begin{itemize}
\item $A=(A_*)'$.
\item $A$ is a Banach algebra, and multiplication in $A$ is
separately $ \sigma (A,A_*)$ continuous; or equivalently $A\cdot A_*\subseteq A_*$ and $A_*\cdot A\subseteq A_*$.
\end{itemize}

In \cite[Theorem 4.1]{J}, it has been proved that for any metric space $(X,d),$
$Lip_\alpha(X,E)$ is a dual space, whenever E is.
In the present paper, we first provide necessary and sufficient conditions for 
that $C_{BSE}(\Delta(A))$ is a Dual Banach algebra. Then we give some conditions,
under which vector-valued Lipschitz algebra $Lip_\alpha(X,\mathcal A)$
is a Dual Banach algebra.

\section{\bf Main Results}

We commence our resuts with the with following propositions.

\begin{pro}
Let $A$ be a commutative Banach algebra such that $\Delta(A)$ be nonempty. Then the following assertions are equivalent. 
\begin{enumerate}
\item[(i)] The weak topology and weak-* topology are the same on $\Delta(A)$, \item[(ii)] $A^{**}|_{\overline{\langle\Delta(A)\rangle}}={\overline{\widehat{A}|_{\overline{\langle\Delta(A)\rangle}}}^{||.||_{A^{**}}}}$, as isometric Banach algebras.
\end{enumerate}
\end{pro}

\begin{proof}
(i)$\Rightarrow$(ii). Since $\widehat{A} \subseteq A^{**}$, it follows that $\widehat{A}|_{\overline{\langle\Delta(A)\rangle}} \subseteq A^{**}|_{\overline{\langle\Delta(A)\rangle}}$ and so 
$\hat{A}|_{\overline{\langle\Delta(A)\rangle}} \subseteq \langle\Delta(A)\rangle^{*}.$
Consequently, $\overline{\widehat{A}|_{\overline{\langle \Delta(A)\rangle}}}^{\|\cdot\|_{A^{**}}} \subseteq {\langle\Delta(A)\rangle}^{*}$. In the other words,
$${\overline{\widehat{A}|_{\overline{\langle\Delta(A)\rangle}}}^{||.||_{A^{**}}}} \subseteq A^{**}|_{\overline{\langle\Delta(A)\rangle}}.$$
Now let $W:= {\overline{\widehat{A}|_{\overline{\langle\Delta(A)\rangle}}}^{||.||_{A^{**}}}}$ and suppose on the contrary that  $W$ is a proper subset of $A^{**}|_{\overline{\langle\Delta(A)\rangle}}$. Thus there exist $F \in A^{**}$  such that  $F|_{\overline{\langle\Delta(A)\rangle}} \notin W$. By Hahn-Banach theorem, there exists $G \in \langle\Delta(A)\rangle^{**}$ such that $G(W)=0$ and $G(F|_{\overline{\langle\Delta(A)\rangle}})=1.$
Let $(x_{\alpha})_{\alpha \in I}$ be a net in $A$ such that $||x_{\alpha}|| \leq ||F||$ $(\alpha\in I)$ and $\hat{x}_{\alpha}\longrightarrow_\alpha F,$ in the weak$^*$ topology of ${\mathcal A}^{**}$. By the assumption, we have $\hat{x}_{\alpha}|_{\langle\Delta(A)\rangle}\longrightarrow_\alpha F|_{\langle\Delta(A)\rangle}$, in the weak topology of ${\mathcal A}^{**}$. Consequently,  $$ 0=G(\hat{x}_{\alpha}|_{\langle\Delta(A)\rangle}) \rightarrow G(F|_{\langle\Delta(A)\rangle}) =1,$$
which is a contradiction.

(ii)$\Rightarrow$ (i). If $\phi_{\alpha}\longrightarrow_\alpha \phi$, in the weak$^*$ topology of $\Delta(A)$ and $F \in A^{**}$, then  $F|_{\overline{\langle\Delta(A)\rangle}} \in {\overline{\widehat{A}|_{\overline{\langle\Delta(A)\rangle}}}^{||.||_{A^{**}}}}$. Thus there exists a sequence $(x_n)_{n\in \mathbb{N}}$ in $A$ such that,  $\hat{{x_{n}}}|_{\overline{<\Delta(A)>}} \longrightarrow_n F|_{\overline{\langle\Delta(A)\rangle}} $, in the norm topology of $A^{**}$.
Consequently for every $\varepsilon>0$ exists $N>0$ such that for all $n \geqslant N$ we have 
$$|| \hat{{x_{n}}}|_{\overline{\langle\Delta(A)\rangle}} - F|_{\overline{\langle\Delta(A)\rangle}} ||<\varepsilon.$$
Thus by taking a subnet if necessary,
$$| \lim_{n} \lim_{\alpha} \phi_{\alpha}(x_{n}) - \lim_{\alpha} F(\phi_{\alpha}) | \leq \varepsilon$$
 and 
$$| \lim_{\alpha} \lim_{n} \phi_{\alpha}(x_{n}) - \lim_{\alpha} F(\phi_{\alpha}) | \leq \varepsilon.$$
It follows that 
$$ \lim_{n} \lim_{\alpha} \phi_{\alpha}(x_{n})=\lim_{\alpha} \lim_{n} \phi_{\alpha}(x_{n})$$
and so $ \lim_{\alpha} F(\phi_{\alpha}) = F(\phi)$. Note that since the nets are bounded, the limits exist.
Consequently, $\phi_{\alpha}\longrightarrow_\alpha\phi,$ in the weak topology of $A^{**}$.
\end{proof}

\begin{pro}
Let $A$ be a commutative Banach algebra such that $\Delta(A)$ be nonempty. Then the following assertions are equivalent. 
\begin{enumerate}
\item[(i)] $A^{**}|_{\overline{\langle\Delta(A)\rangle}}={\overline{\widehat{A}|_{\overline{\langle\Delta(A)\rangle}}}^{||.||_{A^{**}}}}$, as isometric Banach algebras, \item[(ii)] $\overline{A^{**}|_{\langle\Delta(A)\rangle}}^{||.||_{A^{**}}}={\overline{\widehat{A}|_{\langle\Delta(A)\rangle}}^{||.||_{A^{**}}}}$, as isometric Banach algebras.
\end{enumerate}
\end{pro}

\begin{proof}
(i)$\Rightarrow$(ii). It is clear that ${\overline{\hat{A}|_{\langle\Delta(A)\rangle}}^{||.||_{A^{**}}} \subseteq \overline{A^{**}|_{\langle\Delta(A)\rangle}}^{||.||_{A^{**}}}}$. For the proof of the reverse of the inclusion, take $F \in A^{**}|_{\langle\Delta(A)\rangle}$ and define $\overline{F}$ on $\overline{\langle\Delta(A)\rangle}$ as $$\overline{F}(\psi):=\lim_{n}F(\psi_{n}),$$
where, $\{\psi_{n}\}\subseteq\langle\Delta(A)\rangle$ and $\psi_{n}\longrightarrow_n\psi$, in the norm topology.
Note that $\overline{F}$ is well-defined. In fact, for all $m,n\in\Bbb N$ we have
 $$|F(\psi_{n})-F(\psi_{m})| \leq ||F||.||\psi_{n}-\psi_{m}||.$$
 Thus sequence $(F(\psi_n))_{n\in \mathbb{N}}$ is Cauchy and so convergent. Moreover,
 $$|\overline{F}(\psi)|=|\lim_{n}F(\psi_{n})| \leq ||F||.\lim_{n}||\psi_n||=||F||.||\psi||,$$
which implies that $ \overline{F} \in (\overline{<\Delta(A)>})^{*}.$
By the assumption, we have $$A^{**}|_{\overline{<\Delta(A)>}}={\overline{\hat{A}|_{\overline{<\Delta(A)>}}}^{||.||_{A^{**}}}}.$$ It follows that
$\overline{F} \in {\overline{\hat{A}|_{\overline{\langle\Delta(A)\rangle}}}^{||.||_{A^{**}}}}.$ Thus there exists a sequence  $(x_n)_{n\in \mathbb{N}}$ in $A$ such that
$$\hat{{x_{n}}}|_{\overline{\langle\Delta(A)\rangle}}\longrightarrow_n\overline{F},$$
in the norm topology. It follows that
$$\hat{{x_{n}}}|_{\langle\Delta(A)\rangle}\longrightarrow_nF |_{\langle\Delta(A)\rangle},$$ 
in the norm topology.

(ii)$\Rightarrow$(i). Let $A^{**}|_{\langle\Delta(A)\rangle}={\overline{\widehat{A}|_{\langle\Delta(A)\rangle}}^{||.||_{A^{**}}}}$.
It is clear that
$$ {\overline{\widehat{A}|_{\overline{\langle\Delta(A)\rangle}}}^{||.||_{A^{**}}}}  \subseteq A^{**}|_{\overline{\langle\Delta(A)\rangle}}.$$
Now take $H \in  A^{**}|_{\overline{\langle\Delta(A)\rangle}}$ and set
$$ F:= H|_{\langle\Delta(A)\rangle} \in A^{**}|_{\langle\Delta(A)\rangle}.$$
Thus $F \in {\overline{\widehat{A}|_{\langle\Delta(A)\rangle}}^{||.||_{A^{**}}}}$ and so there exists a sequence $(x_n)_{n\in \mathbb{N}}$ in $A$ such that $\hat{{x_{n}}}|_{\langle\Delta(A)\rangle}\longrightarrow_n F$,
in the norm topology. It follows that $\hat{{x_{n}}} (\psi)\longrightarrow_nF(\psi)$,  for all $\psi \in \langle\Delta(A)\rangle$ with $\|\psi\|\leq 1$.
Now take $k \in \overline{\langle\Delta(A)\rangle}$ such that $\|k\|\leq 1$. There exists the sequence $(k_n)_{n\in \mathbb{N}}$ in $\langle\Delta(A)\rangle$ such that, $k_{n}\longrightarrow_nk$, in the norm topology. Thus there exists $N\in\Bbb N$ such that for all $n,\ell \geq N$, $\|k_{\ell}\|\leq 1$ and
$$|\hat{{x_{n}}} (k_{\ell}) -F(k_{\ell})|< \varepsilon.$$
Consequently,
$$|\hat{{x_{n}}} (k)-F(k)|=|\lim_{\ell}(\hat{{x_{n}}} (k_{\ell}) -F(k_{\ell}))|=\lim_{\ell}|\hat{{x_{n}}} (k_{\ell}) -F(k_{\ell})| \leq \varepsilon.$$
It follows that for all $n \geq N,$
$$||\hat{{x_{n}}}|_{\overline{\langle\Delta(A)\rangle}}-F ||<\varepsilon,$$
which implies that $F \in {\overline{\widehat{A}|_{\overline{\langle\Delta(A)\rangle}}}^{||.||_{A^{**}}}}.$ 
Therefore,
$$A^{**}|_{\overline{\langle\Delta(A)\rangle}}  \subseteq {\overline{\widehat{A}|_{\overline{\langle\Delta(A)\rangle}}}^{||.||_{A^{**}}}},$$
and the proof is completed.
\end{proof}

We state here the main result of the present paper.

\begin{thm}\label{t1}
Let $A$ be a commutative Banach algebra such that $\Delta(A)$ be nonempty. Then the following assertions are equivalent. 
\begin{enumerate}
\item[(i)]  $C_{BSE}(\Delta(A))=\langle\Delta(A)\rangle^*$ ,as isometric Banach algebras, \item[(ii)] The weak topology and weak-* topology are the same on $\Delta(A)$.
\end{enumerate}
\end{thm}

\begin{proof}
(i)$\Rightarrow$(ii). Suppose that $\{\phi_{\alpha}\}$ is a net in $\Delta(\mathcal A)$, converging to $\phi\in \Delta(A)$, in the weak$^*$ topology of $\Delta(\mathcal A)$.  If  $F \in A^{**}$, then $F|_{\overline{\langle\Delta(A)\rangle}} \in \langle\Delta(A)\rangle^*$ and so $F|_{\overline{\langle\Delta(A)\rangle}} \in C_{BSE}(\Delta(A))$.
It follows that $F|_{\Delta(A)} \in C_{BSE}(\Delta(A))$ and we have $$|| F|_{\Delta(A)} ||_{BSE} = || F|_{\langle\Delta(A)\rangle} ||_{{A}^{**}} \leq ||F|| < \infty.$$
Since $F|_{\Delta(A)} \in C_{b}(\Delta(A))$, thus $$F|_{\Delta(A)}(\phi_{\alpha}) \longrightarrow_\alpha F|_{\Delta(A)}(\phi).$$ 
Consequently $F(\phi_{\alpha}) \rightarrow_\alpha F(\phi)$, which implies that $\{\phi_{\alpha}\}$ converges to  $\phi$,
in the weak topology of $\Delta(\mathcal A)$.
(ii)$\Rightarrow$ (i). Define $$\Phi : C_{BSE}(\Delta(A)) \longrightarrow\langle\Delta(A)\rangle^{*},$$ by 
$$\Phi (H)( \sum_{k=1}^{n} \lambda_{k} \phi_{k}):= \sum_{k=1}^{n} \lambda_{k} H(\phi_{k})\;\;\;\;\;(H \in  C_{BSE}(\Delta(A))).$$ It is obvious that $\Phi$ is linear and
$$\left|\Phi (H)( \sum_{k=1}^{n} \lambda_{k} \phi_{k})\right|\leq ||H||_{BSE} \;|| \sum_{k=1}^{n} \lambda_{k} \phi_{k} ||_{A^{*}}.$$ Thus $||\Phi(H)|| \leq  ||H||_{BSE}$ and so $\Phi$ is continuous.
For any $G \in \langle\Delta(A)\rangle^{*}$ there exists $\sigma:=G|_{\Delta(A)}$ such tha  $\Phi(\sigma)=G.$ In fact,
\begin{equation}\label{e1}
\Phi (\sigma)( \sum_{k=1}^{n} \lambda_{k} \phi_{k})=\sum_{k=1}^{n} \lambda_{k} G(\phi_{k})=G(\sum_{k=1}^{n} \lambda_{k} \phi_{k}),
\end{equation}
which implies that $\Phi(\sigma)=G$. Now we show that $\sigma$ is. Let $\phi_{\alpha} \longrightarrow_\alpha\phi$, in the weak$^*$ topology of $\Delta(A)$. By the assumption, $\phi_{\alpha} \longrightarrow_\alpha\phi$,
in the weak topology of $\Delta(A)$ and so $G(\phi_{\alpha}) \rightarrow G(\phi) $. Thus $\sigma(\phi_{\alpha}) \rightarrow \sigma(\phi) $. Consequently $\sigma \in C_{b}(\Delta(A))$ and $$|\sum_{k=1}^{n}c_{i} \sigma(\phi_{i})| \leq ||G||\;||\sum_{k=1}^{n}c_{i} \phi_{i}||_{A^{*}}$$
which implies  that $\sigma \in C_{BSE}(\Delta(A))$ and $||\sigma||_{BSE} \leq ||G||.$ Thus $\Phi$ is surjective.
Moreover, by \cite[]{}, for each $H \in  C_{BSE}(\Delta(A))$, $H=F|_{\Delta(A)}$ for some $F\in A^{**}$. Then
\begin{eqnarray*}
||H||_{BSE}&=&\sup \left\{ |\sum_{k=1}^{n}\lambda_{k}H(\phi_{k})| : ||\sum_{k=1}^{n} \lambda_{k} \phi_{k}||_{A^{*}} \leq 1 \right\}\\
&=&\sup \left\{ |F(\sum_{k=1}^{n}\lambda_{k}\phi_{k})|||\sum_{k=1}^{n} \lambda_{k} \phi_{k}||_{A^{*}} \leq 1 \right\}\\
&=&\sup \left\{ |F|_{<\Delta(A)>}(\psi)|:||\psi||_{<\Delta(A)>} \leq 1 \right\}\\
&=&||F|_{<\Delta(A)>}||_{A^{**}}\\
&=&||\Phi(H)||_{A^{**}}.
\end{eqnarray*}
Therefore $\Phi$ is isometric. Now we show that $\Phi$ is isomorphism. For $\sigma_{1} , \sigma_{2} \in  C_{BSE}(\Delta(A))$, there exist $F_{1} , F_{2} \in A^{**}$ such that $\sigma_{1}= F_{1}|_{\Delta(A)} , \sigma_{2}= F_{2}|_{\Delta(A)}$. Then for each $\phi \in \Delta(A)$, we have
$$(\sigma_{1}.\sigma_{2})(\phi):=\sigma_{1}(\phi).\sigma_{2}(\phi)=(F_{1}|_{\Delta(A)})(\phi).(F_{2}|_{\Delta(A)})(\phi)=F_{1}(\phi)F_{2}(\phi).$$
Furthermore,
$$[\Phi(\sigma_{1}) \Box \Phi(\sigma_{2})](\phi)=\Phi(\sigma_{1})[\Phi(\sigma_{2})\phi]$$ 
and so for each $x \in A,$
$$\Phi(\sigma_{2})\phi(x)=\Phi(\sigma_{2})(\phi.x),$$
where, $$(\phi.x)(y)=\phi(xy)=\phi(x)\phi(y)\;\;\;\;\;(y\in A).$$
Thus for all $x \in A $ and $ \phi \in \Delta(A)$, $ \phi.x=\phi(x)\phi.$
Consequently,
$$\Phi(\sigma_{2})\phi(x)=\Phi(\sigma_{2})(\phi(x)\phi) =\phi(x)\Phi(\sigma_{2})(\phi)$$
and so for any $\sigma_{2} \in C_{BSE}(\Delta(A))$,
$$\Phi(\sigma_{2})\phi=\Phi(\sigma_{2})(\phi)\phi.$$
It follows that
$$[\Phi(\sigma_{1}) \Box \Phi(\sigma_{2})](\phi)=\Phi(\sigma_{1})[\Phi(\sigma_{2})(\phi)\phi]=\Phi(\sigma_{2})(\phi)\Phi(\sigma_{1})(\phi).$$
It follows that
$$[\Phi(\sigma_{1}) \Box \Phi(\sigma_{2})](\phi)=F_{1}(\phi)F_{2}(\phi)=\sigma_{1}(\phi)\sigma_{2}(\phi)=\Phi(\sigma_{1}\sigma_{2})(\phi).$$
Therefore
$$\Phi(\sigma_{1}\sigma_{2})=\Phi(\sigma_{1}) \Box \Phi(\sigma_{2}).$$
and so $\Phi$ is isomorphism.
\end{proof}

\begin{cor}
Let $A$ be a reflexive commutative Banach algebra such that $\Delta(A)$ be nonempty. Then 
 $C_{BSE}(\Delta(A))=\langle\Delta(A)\rangle^*$.
\end{cor}

\begin{exa}\em
Let $(X,d)$ be an infinite compact metric space and $A=Lip_dX$. Then 
$$
\Delta(A)=\left\{\delta_x : x \in X \right\}.
$$
Moreover, we have
$$
A=\langle\delta_x : x \in X \rangle^*.
$$
Since $A$ is a unital BSE-algebra \cite[Examp,e 6.1]{Kan}, we obtain
$$C_{BSE} (\Delta (A))=\widehat{M(A)}=\widehat{A}=\langle\delta_x: x \in X\rangle^*=\langle\Delta(A)\rangle^*.$$
Note that $A$ is not reflexive. In fact, by \cite{D},
$$lip_d (X)^{**}=Lip_d X\;\;\;\text{and}\;\;\;\Delta (lip_d X)=\Delta (Lip_d X).$$
Thus $A$ is reflexive if and only if $lip_d (X)^{**} $ is reflexive. It follows that 
${lip_d X}^*$ and so $lip_d X $ is reflexive. Consequently, $Lip_d X=lip_d X$,
which implies that  $X$ is descrete uniform. Since $X$ is compact it follows that $X$ is
finite, which is a contradiction. 
\end{exa}

\begin{pro}
Let $A$ be a commutative and semisimple Banach algebra such that 
weak topology and weak$^*$ topology be the same on $ \Delta(A)$. Then
$C_{BSE}(\Delta(A)) $ is a dual Banach 
algebra by predual $B={\overline{\langle\Delta(A)\rangle}}^{||.||_{A^*}}$.
\end{pro}

\begin{proof}
According to Theorem \ref{t1},
$$C_{BSE}(\Delta(A)) \cong \langle\Delta(A)\rangle^*={\overline{\langle\Delta(A)\rangle}}^*= B^*,$$
as two isometric Banach algebras. It is sufficinet to show that
$$C_{BSE}(\Delta(A))\cdot\Delta(A)\subseteq B.$$
Let $\sigma \in C_{BSE}(\Delta(A))$ and $ \phi \in \Delta(A)$. Then there exists $F \in A^{**}$ such that
$\sigma = F|_{\Delta(A)}$. Then
$$(\sigma\cdot\phi )(x)=\sigma(\phi\cdot x)=F( \phi\cdot x)\;\;\;\;\;\;(x\in A).$$
Thus
$$ \phi\cdot x=\phi (x)\phi,$$
and so
$$(\sigma\cdot\phi )(x)=F(\phi (x)\phi)=\phi(x)F( \phi ).$$
It follows that
$$\sigma\cdot\phi=F( \phi )\phi\in\langle\Delta(A)\rangle\subseteq B.$$
In general, if $\psi=\Sigma_{i=1}^{n} \lambda_i \phi_i \in\langle\Delta(A)\rangle$, where $\phi_i  \in \Delta(A),$ 
$(i=1,\cdots,n)$, we have
$$\sigma\cdot\psi=\sum_{i=1}^{n} \lambda_i F(\phi_i) \phi_i \in B.$$
Now, take $\psi \in B.$ There exists a sequence $\{\psi_n\}\in\langle\Delta(A)\rangle$ such that, 
$\psi_n \longrightarrow_n\psi $, in the norm topology. Thus
$$|| \sigma\cdot\psi - \sigma\cdot\psi_n|| \leq ||\sigma ||||\psi - \psi_n ||\rightarrow_n0,$$
and so $\sigma\cdot\psi\in B.$ Therefore  $ (C_{BSE}(\Delta(A)) , B)$ is a dual Banach algebra.
\end{proof}

\begin{thm}\label{t2}
Let $(X,d)$ be a metric space and $E$ be a Banach algebra. Suppose that 
$$W=\overline{\langle\delta_s \otimes f: \;f \in E^* , s \in X\rangle} \subseteq Lip_d(X,E)^*,$$ 
where,
$$(\delta_s \otimes f)(h)=f(h(s))\;\;\;\;\;\;\;(s \in X,f \in E^*,h \in Lip_d(X,E)).$$
Then $Lip_d(X,E) $ can be imbedded isomerically in the dual Banach algebra $W^{*}$.
\end{thm}

\begin{proof}
Define
$$\Phi : Lip_d(X,E) \rightarrow W^{*}$$
$$h \rightarrowtail \Phi_{h},$$
where
$$\Phi_{h}( \Sigma_{k=1}^{n} \lambda_k (\delta_{s_k} \otimes f_k)):=\Sigma_{k=1}^{n} \lambda_k f_k(h(s_k)),$$
for each $h \in Lip_d(X,E).$ It is clear that $\Phi $ is well-defined and linear. 
We show that $\Phi$ is an isometric mapping. It is easily verified that $||\Phi_h || \leq ||h||.$ For the reverse inequality,
note that for each $h\in Lip_d(X,E)$ and $x\in X$ we have
\begin{eqnarray*}
||h(x)||_E&=&||\widehat{h(x)}||_{E^{**}}=\sup \left\{|\sigma(h(x))|: \sigma \in E^{*}, ||\sigma|| \leq 1 \right\}\\
&=&\sup\left\{ | (\delta_x\otimes \sigma )(h) |: \sigma \in E^{*}, ||\sigma|| \leq 1 \right\}.
\end{eqnarray*}
Moreover, for any $x \in X$ and $ \sigma \in E^*$ with $||\sigma|| \leq 1$ we have
$||\delta_x\otimes \sigma||\leq 1$ and 

\begin{eqnarray*}
||\Phi_h ||&=&\sup \left\{|\Phi_{h}(u)|: \;u \in W, ||u||\leq 1\right\}\\
&=&\sup\left\{\left|\Phi_{h}(\sum_{k=1}^{n}\delta_{x_k}\otimes f_k)\right|: 
\left\|\sum_{k=1}^{n}\delta_{x_k}\otimes f_k\right\|\leq 1\right\}\\
&\geq &|\Phi_h(\delta_{x}\otimes \sigma)|.
\end{eqnarray*}

Thus we obtain

\begin{eqnarray*}
||h(x)||_E &=&\sup \left\{|\sigma(h(x))|: \sigma \in E^{*}, ||\sigma|| \leq 1 \right\}\\
 &=&\sup \left\{|\delta_x\otimes \sigma (h) | : \sigma \in E^{*} , ||\sigma|| \leq 1 \right\}\\
&=&\sup \left\{|\Phi_h ( x \otimes \sigma ) : \sigma \in E^{*} , ||\sigma|| \leq 1 \right\}\\
&\leq&\sup\left\{|\Phi_h ( \sum_{k=1}^{n}\delta_{x_k}\otimes f_k)|: ||\sum_{k=1}^{n}
\delta_{x_k}\otimes f_k|| \leq 1 \right\}\\ &=&||\Phi _h||.
\end{eqnarray*}

 Consequently,
 $$||h||_{\infty , E} \leq || \Phi _h ||.$$
Moreover, for all $x,y\in X$ with $x\neq y$ we have
\begin{eqnarray*}
\frac{||h(x)-h(y)||_E }{ d(x,y) }&\leq&\frac{||h(x)-h(y))||_E }{||\delta_x -\delta_y || }\\
&=&\sup \left\{ \frac{|\sigma(h(x))-\sigma(h(y)) | }{||\delta_x -\delta_y ||}: \sigma\in E^*,\|\sigma\|\leq 1\right\}\\
&=&\sup \left\{ \frac{||\delta_{x}\otimes \sigma (h)- \delta_{y}\otimes \sigma (h) || }{|| \delta_x -\delta_y ||}:  \sigma\in E^*,
\|\sigma\|\leq 1\right\}.
\end{eqnarray*} 
On the other hand, for any $\sigma\in E^*$ with $\|\sigma\|\leq 1$
\begin{eqnarray*}
||\delta_x \otimes \sigma -\delta_y\otimes \sigma||&=&\sup \left\{| (\delta_x \otimes \sigma -\delta_y\otimes \sigma) (h) | : h \in Lip_d(X,E) , ||h|| \leq 1 \right\}\\
 &=&\sup\left\{ | \sigma(h(x))- \sigma(h(y)) | :  h \in Lip_d(X,E) , ||h|| \leq 1 \right\}\\
&\leq& ||h(x) - h(y) ||\\
&=& ||(\delta_x -\delta_y)(h)||\\
&\leq& || \delta_x -\delta_y ||||h ||\\
&\leq& || \delta_x -\delta_y ||.
\end{eqnarray*} 
Consequently,
$$\frac{1}{|| \delta_x -\delta_y ||} \leq  \frac{1}{||\delta_x\otimes \sigma -\delta_y\otimes \sigma ||}.$$
Thus
\begin{eqnarray*}
\frac{||h(x)-h(y)||_E}{d(x,y)}&\leq&\sup \left\{ \frac{|\delta_x\otimes \sigma(h)-\delta_y\otimes \sigma (h) | }{|| \delta_x -\delta_y ||}:  ||\sigma ||\leq 1\right\}\\
&\leq&\sup \left\{ \frac{| (x \otimes \sigma - y \otimes \sigma) (h) | }{||\delta_x\otimes \sigma - y \otimes \sigma ||}: || \sigma ||\leq 1 \right\}\\
&\leq&\sup\left\{ \frac{ ||\Phi_h (u) ||}{||u||} : u \neq 0 \right\}\\
&=&||\Phi_h ||.
\end{eqnarray*} 
Then for all $ h \in Lip_d(X,E)$, we  obtain
$$\rho_{d,E} (h)=\sup \left\{\frac{||h(x)- h(y) ||_E }{d (x, y)}:  x \neq y\right\}\leq ||\Phi_h ||.$$
Then
$$|| h ||=\max \left\{\rho_{d,E} (h), ||h||_{\infty , E}\right\}\leq||\Phi_h|| .$$
Therefore $\Phi$ is isometry, as claimed.
\end{proof}

\begin{pro}\label{p2}
Let $(X,d)$ be a metric space and $E$ be a natural Banach algebra. Then for each $f\in Lip_d(X,E)$ we have
$$||f||=||\widehat{f}||=\sup\left\{|h(f) |: h \in E^*, ||h|| \leq 1\right\}=\sup \left\{ |\phi(f) |: \phi\in\Delta (E) \right\}.$$
\end{pro}

\begin{proof}
Since $E$ is natural, $\Delta (E)=\left\{\delta_t : t \in X \right\}$.
It is clear that,
$$\sup \left\{|\phi(f)|: \phi \in \Delta (E) \right\}=\sup \left\{|\delta_t (f)|: t \in X \right\}
=||f||_\infty\leq ||f||.$$
Moreover, 
$$|\Phi_{f}(u) | = | \Sigma_{k=1}^{n}  \phi_k(f(x_k)) | \leq ||u|| . ||f||_{Lip_d(X,E)}.$$
Thus
$$||\Phi_f|| \leq  ||f||_{Lip_d(X,E)}.$$ 
Since $u =\Sigma_{k=1}^{n}  \delta_{x_k} \otimes \phi_k  \in W \subseteq (Lip_d(X,E))^{*},$ thus
$$||u|| =\sup \left\{ | \Sigma_{k=1}^{n}  ( \delta_{x_k} \otimes \phi_k )(f) | : f \in Lip_d(X,E) , ||f|| \leq 1 \right\} $$
$$=\sup \left\{ | \Sigma_{k=1}^{n} \phi_k(f_{x_k})) | : f \in Lip_d(X,E) , ||f|| \leq 1 \right\} .$$
Now we show that $P_{d,E}(f) \leq ||\Phi_f ||$ and $||f||_{\infty , E} \leq ||\Phi_f ||.$
Note that $$\|f||_{\infty , E}=\sup \left\{ ||f(s)||_E : s \in X \right\},$$ where 
$$||f(s)||_E=|| \hat{f(s)}||_{E^{**}}=\sup \left\{ |h(f(s))| : h \in E^{*}, ||h|| \leq 1 \right\}$$
and
$$||\Phi_f||=\sup\left\{\left|\sum_{k=1}^{n} \phi_k(f_{x_k}))\right|: \left\| \sum_{k=1}^{n} \delta_{x_k} \otimes \phi_k \right\| \leq 1 \right\}.$$
By the hypothesis, for any $z \in E,$ 
$$||z||=\sup\left\{ | \phi(z) | : \phi \in \Delta(E) \right\}.$$
Thus
\begin{eqnarray*}
||f(s)||&=&\sup\left\{ | \phi(f(s)) | : \phi \in \Delta(E) \right\}\\
&=&\sup \left\{ |(\delta_{s} \otimes \phi )(f) | : \phi \in \Delta(E) \right\}.
\end{eqnarray*}
Since
$$|| \delta_{s} \otimes \phi||= ||\delta_{s} ||\;|| \phi||=1$$
and for every $s \in X$ 
$$||f(s)|| \leq ||\phi(f)||,$$
it follows that
$$||f||_{\infty , E} \leq ||\Phi_f ||.$$
Furthermore, for all $s,t\in X$ with $s\neq t$ we have
\begin{eqnarray*}
\frac{||f(s)-f(t)||_E }{d(s,t) }&\leq&\frac{||f(s)-f(t)||_E }{||\delta_s -\delta_t || }\\
&=&\sup \left\{ \frac{| \phi (f(s))- \phi (f(t)) | }{ ||\delta_s -\delta_t || }: s \neq t \right\}\\
&=&\sup \left\{ \frac{||( \delta_s \otimes \phi - \delta_t \otimes \phi)(f) || }{ || \delta_s -\delta_t || }: s \neq t \right\}\\ &\leq&\sup \left\{ \frac{|\Phi_f(u)|}{||u||} : u \neq 0 \right\}\\
&=&||\Phi_f ||.
\end{eqnarray*}
Therefore $\rho_{d,E}(f) \leq ||\Phi_f ||$ and so $\Phi $ is isometry.
\end{proof}

\begin{pro}\label{p4}
Let $Y$ be a Hausdorff locally compact space, $(X,d)$ be a metric space and $E=(C_0(Y) , ||.||_\infty).$ 
Suppose that 
$$W= \overline{\langle\delta_x \otimes \delta_y: x \in X , y \in Y\rangle}.$$
Then 
$$Lip_d(X,E)=\langle\delta_x\otimes \delta_y: x\in X , y \in Y\rangle^{*}:= W^{*} $$
as a dual Banach algebra. 
\end{pro}

\begin{proof}
Define 
$$\Phi : Lip_d(X,E) \rightarrow W^*$$
$$f \rightarrowtail \Phi_{f},$$
where
$$ \Phi_{f}( \Sigma_{k=1}^{n} \lambda_k (\delta_{x_k} \otimes \delta_y))=\sum_{k=1}^{n} \lambda_k f(x_k) (y_k).$$
It is easily verified that $\Phi$ is linear, and $$|| \Phi_f || \leq ||f||.$$
We prove the reverse of the inequality. For all $s,t \in X$ with $s\neq t$ we have
\begin{eqnarray*}
\frac{||f(s)-f(t)||_E }{ d(s,t) }&\leq&\frac{||f(s)-f(t)||_\infty }{|| \delta_s -\delta_t ||_{Lip_d(X,E)^{*}} }\\
&=&\sup \left\{ \frac{| f(s)(y)- f(t)(y) |}{||\delta_s -\delta_t ||_{Lip_d(X,E)^{*}} }: y \in Y \right\}\\
&\leq&\sup \left\{ \frac{|\Phi_f(u)|}{||u||}: u \in W , u \neq 0 \right\}.
\end{eqnarray*}
Since
\begin{eqnarray*}
|f(s)(y)- f(t)(y)|&=&|\Phi_f(\delta_s \otimes \delta_y)-\Phi_f(\delta_t \otimes \delta_y) |\\
&=& |\Phi_f((\delta_s-\delta_t ) \otimes \delta_y) |\\
&\leq&||\Phi_f||\;||(\delta_s - \delta_t ) \otimes \delta_y ||\\
&\leq& ||\Phi_f||\;|| \delta_s - \delta_t ||||\delta_y ||\\
&=&||\Phi_f||\;|| \delta_s - \delta_t ||,
\end{eqnarray*}
it follows that
\begin{equation}\label{e2}
\frac{| f(s)(y)- f(t)(y) | }{ || \delta_s -\delta_t || } \leq ||\Phi_f||.
\end{equation}
For all $ s,t \in X$ with $s\neq t$ we have
\begin{eqnarray*}
|| \delta_{s} - \delta_{t} ||&=&\sup \left\{|\delta_{s}(f)-\delta_{t}(f) | : f \in Lip_d(X,E), |||f|||_{d,E} \leq 1 \right\}\\
&=&\sup\left\{|f(s)-f(t)|: f \in Lip_d(X,E) , |||f|||_{d,E} \leq 1 \right\}\\
&\leq&\sup\left\{P_{d,E}(f)\;d(s,t) : f \in Lip_d(X,E), |||f|||_{d,E} \leq 1 \right\}\\
&\leq& d(s,t).
\end{eqnarray*}
This inequality together with \eqref{e2} imply that
$$\frac{||f(s)-f(t)||_E }{ d(s,t) } \leq ||\Phi_f||.$$
Consequently,
$$\rho_{d,E}(f) \leq ||\Phi_f||=||\Phi(f)||.$$
Moreover, it is easily verified that
$$||f||_{\infty , E}\leq ||\Phi_f||.$$
It follows that
$$\max\left\{ ||f||_{\infty , E},\rho_{d,E}(f)  \right\}=||f||\leq ||\Phi(f)||,$$
and so $\Phi $ is isometry.  Now we show that $\Phi $ is surjective. 
Take $F \in W^*$ and for every $x\in X$ and $y \in Y$, define
$$f(x)(y):= F(\delta_{x} \otimes \delta_{y}).$$
Then for all $ s,t\in X$ with $s\neq t$ we have
$$ \frac{||f(s)-f(t)||_E }{d(s,t)} \leq ||F||$$
and so 
$$\rho_{d,E}(f) \leq  ||F||.$$
Also
$$||f||_{\infty,E} \leq ||F||,$$
obviously. Consequently,
$||f\|\leq||F||$ and so $f\in  Lip_d(X,E).$
It is not hard to see that $\Phi(f)=F$. Thus $\Phi$ is surjective. 
Finally, we show that 
$$Lip_d(X,E)\cdot W \subseteq W\;\;\text{and}\;\;W\cdot Lip_d(X,E) \subseteq W.$$
Let $f \in Lip_d(X,E) $, $y,z \in Y$ and $s,t \in X$. For every $g \in  Lip_d(X,E) $ we have
$$ ((\delta_{s}\otimes\delta_{y})\cdot\Phi_f )(g) = \Phi_f(g\cdot(\delta_{s}\otimes \delta_{y})). $$
But
\begin{eqnarray*}
g\cdot(\delta_{s}\otimes\delta_{y})(f)&=&(\delta_{s}\otimes \delta_{y})(fg)\\
&=&(fg(s))(y)\\
&=&(f(s)g(s))(y)\\
&=&f(s)(y)\;g(s)(y).
\end{eqnarray*}
Thus
$$g\cdot (\delta_{s} \otimes\delta_{y})(f)=(\delta_{s} \otimes \delta_{y})(f)\;(\delta_{s}\otimes \delta_{y})(g),$$
and so
$$g\cdot (\delta_{s} \otimes \delta_{y})=(\delta_{s} \otimes \delta_{y})(g)\cdot (\delta_{s} \otimes \delta_{y}).$$
Thus
\begin{eqnarray*}
[(\delta_{s} \otimes \delta_{y})\cdot\Phi_f](g)&=&\Phi_f[(\delta_{s} \otimes \delta_{y})(g)\cdot(\delta_{s}\otimes \delta_{y})]\\
&=&(\delta_{s} \otimes \delta_{y})(g)\cdot\Phi_f(\delta_{s}\otimes\delta_{y}).
\end{eqnarray*}
Consequently,
$$(\delta_{s} \otimes \delta_{y})\cdot\Phi_f =\Phi_f(\delta_{s} \otimes \delta_{y})\cdot(\delta_{s} \otimes \delta_{y})$$
and so
$$(\delta_{s} \otimes \delta_{y})\cdot\Phi_f = f(s)(y).(\delta_{s} \otimes \delta_{y}) \in W.$$
Now let  $u=\Sigma_{k=1}^{n} \lambda_k \delta_{s_k} \otimes \delta_{y_k} \in W $. Then
$$u\cdot\Phi_f=\sum_{k=1}^{n} \lambda_k f(s_k)(y_k)(\delta_{s_k} \otimes \delta_{y_k}) \in W.$$
In general, if $ u_n \rightarrow u $, in the norm topology, where $u_n\in\langle\delta_{s} \otimes \delta_{y}\rangle,$ then
$$||u_n\cdot\Phi_f - u\cdot\Phi_f|| \leq ||u_n - u ||\;||\Phi_f|| \rightarrow 0.$$
Therefore
$$W\cdot Lip_d(X,E)\subseteq W.$$
Similarly, one can show
$$Lip_d(X,E)\cdot W\subseteq W,$$
and so the proof is completed.
\end{proof}

\begin{pro}\label{p1}
Let $(X,d)$ be a metric space and $E$ be a dual Banach algebra such that $E= F^*$. Then
$$Lip_d(X,E) \cong \langle X \otimes F\rangle^{*}.$$
\end{pro}

\begin{proof}
Define
$$\Phi : Lip_d(X,E) \rightarrow V^*$$
by $f\mapsto\Phi_f$, where
$$V=\langle x \otimes z : x \in X , z \in F\rangle,$$
and for each $f \in Lip_d(X,E),$
$$\Phi_f(x \otimes z) = f(x)(z) . $$
Johnson \cite{J} showed that the mapping of $\Phi $ is linear and isometric and
$$ Lip_d (X,F^*) = V^*,$$
as two Banach spaces. Moreover, for all $f, g \in Lip_d(X,E)$ we have
$$\Phi_{f\;g}(x \otimes z)=(f\;g)(x)(z)=f(x)(z)\;g(x)(z)=\Phi_{f}(x \otimes z)\;\Phi_{g}(x \otimes z) .$$
Thus $\Phi $ is a homomorphism. Finally, we show that
$$ Lip_d(X,E)\cdot V \subseteq V\;\;\text{and}\;\; V\cdot Lip_d(X,E)  \subseteq V.$$
Let $ k \in Lip_d (X,E)$ and $ x \otimes z \in V$. Then for each $h\in Lip_d (X,E)$ we have
\begin{eqnarray*}
 ((x \otimes z).k)(h)&=&(x \otimes z) (kh)=(kh(x))(z)\\
&=&(k(x)h(x))(z)=k(x)(z)h(x)(z)\\
&=&(k(x)(z))((x\otimes z) (h)).
\end{eqnarray*}
It follows that
$$ (x \otimes z)\cdot k=(k(x)(z)) (x \otimes z).$$
Thus
\begin{eqnarray*}
\Phi_{f} ((x \otimes z)\cdot k)&=&\Phi_{f}( (k(x)(z)) (x \otimes z))\\
&=&k(x)(z) \Phi_{f} (x \otimes z)\\
&=&(x \otimes z)(k) \Phi_{f} (x \otimes z)\\
&=&(\Phi_f\cdot(x\otimes z) )(k),
\end{eqnarray*}
which implies that
$$\Phi_f\cdot(x \otimes z)=(\Phi_f (x \otimes z))\cdot(x \otimes z) \in V.$$
Similarly, one can prove that 
$$(x \otimes z)\cdot\Phi_f=(x \otimes z) \cdot (\Phi_f (x \otimes z))\in V,$$
and so the proof is completed.
\end{proof}

\begin{rem}\rm
Proposition \ref{p1} asserts that if $E$ is a dual Banach algebra, 
then $Lip_d(X,E)$ is a dual Banach algebra, as well. 
But the converse of this statement is not true, by Proposition \ref{p4}.
\end{rem}

\vspace{9mm}

{\footnotesize \noindent
F. Abtahi\\
Faculty of Mathematics and Statistics\\
Department of Pure Mathematis\\
University of Isfahan\\
Isfahan, Iran\\
f.abtahi@sci.ui.ac.ir\\
abtahif2002@yahoo.com\\

\noindent
 A. Rejali\\
Faculty of Mathematics and Statistics\\
Department of Pure Mathematis\\
University of Isfahan\\
Isfahan, Iran\\
rejali@sci.ui.ac.ir\\

\noindent
 F. Sayaf\\
Faculty of Mathematics and Statistics\\
Department of Pure Mathematis\\
University of Isfahan\\
Isfahan, Iran\\
f.sayaf@sci.ui.ac.ir\\

\end{document}